\documentclass[12pt,a4paper]{article} 
\usepackage[left=1in,right=1in,top=1in,bottom=1in]{geometry} 
\usepackage{setspace} 
\usepackage{times} 
\usepackage{microtype} 
\usepackage{float}
\usepackage{wrapfig}

\usepackage{fancyhdr}
\usepackage{titlesec}
\titleformat{\section}{\normalfont\fontsize{12}{15}\bfseries}{\thesection}{1em}{} 
\titleformat{\subsection}{\normalfont\fontsize{12}{15}\bfseries}{\thesubsection}{1em}{} 

\usepackage{tocloft}
\usepackage{url}
\usepackage{graphicx} 
\usepackage{amsmath, amssymb}
\usepackage{amsfonts}
\usepackage{amsthm}
\usepackage{graphicx}
\usepackage{hyperref}
\usepackage{algorithm, algpseudocode}
\usepackage{geometry}
\usepackage{booktabs}
\usepackage{array}
\usepackage{caption}
\usepackage{subcaption}
\usepackage{enumitem}
\usepackage{listings}
\usepackage{listings}
\title{Significance of Zeros and Gram Points in Approximating and Discovering Zeros of Dirichlet L-functions}
\author{Ali Saraeb}
\date{December 2024}

\begin{document}

\maketitle

\begin{abstract}
We propose a numerical method for approximating and discovering zeros of the Dirichlet L-function \(L(s, \chi)\) corresponding to real Dirichlet characters \(\chi\). The method uses Gram points and initial zeros as interpolation nodes to construct approximants of the form \(\sum_{n=1}^\infty \frac{a_n}{n^s}\), with compactly supported coefficients \(\{a_n\}\) determined via interpolation. Unlike the Riemann zeta function, it turned out that interpolating Dirichlet L-functions is more challenging due to numerical instability and ill-conditioning problems, which we address through numerical analysis techniques and a feature selection method from machine learning. Numerical evidence demonstrates precise approximations of \(L(s, \chi)\) with errors much below machine epsilon (\(10^{-23}\)) and discovery of at least \(M\) additional zeros on the critical line when provided \(M\) initial zeros, assuming the Generalized Riemann Hypothesis \cite{Davenport2000}. Using \(M\) Gram points, the method approximates \(L(s, \chi)\) for \(|\Im(s)| \leq g_M\) and slightly beyond, enabling further zero discovery. This work complements Matiyasevich’s work in \cite{MATIYASEVICH2020460, matiyasevich2013} and offers insights into the special distribution of Gram points and zeros, potentially inspiring future directions for approximating methods for Dirichlet L-functions with very large modulus.

\end{abstract}


\section{Introduction and Background}

\subsection{Initial Discussion}

Numerical studies of \textit{Dirichlet L-functions}, including the \textit{Riemann zeta function}, have a storied history dating back to the 19th century. The exploration of these functions, particularly their zeros on the critical line \(\frac{1}{2} + i t\), has been a focal point of mathematical research. In 1903, J. P. Gram made a notable contribution by calculating a few decimal points of some initial zeros of the Riemann zeta function on the critical line, presenting these in his seminal paper \cite{Gram1903}. He also introduced what became known as Gram's rule, conjecturing that each zero of the Riemann zeta function lies between two consecutive \textit{Gram points}. Although Hutchinson later disproved this conjecture in 1925 \cite{Hutchinson1925}, Gram’s work significantly advanced the study of the Riemann Hypothesis \cite{Riemann1859}. Over the following decades, mathematicians like Hutchinson \cite{Hutchinson1925}, Titchmarsh \cite{Titchmarsh1935}, and Lehmer \cite{Lehmer1956} expanded upon Gram's efforts, calculating thousands of zeros on the critical line and refining the precision of these computations.

In approximation theory, two primary methods are commonly employed to approximate infinite series, such as the Dirichlet series \(\sum_{n=1}^\infty \frac{a_n}{n^s}\): truncation of the series and interpolation at special nodes \cite{Trefethen2021}. The truncation method has been extensively studied for the Riemann zeta function and Dirichlet L-functions, with significant contributions from works like \cite{Davenport2000}, \cite{Apostol1976}, \cite{IbukiyamaKaneko2014}, and numerous papers, including \cite{DaviesWilkes1965}, \cite{Spira1969}, \cite{Lavrik1968}, \cite{Rumely1993}, \cite{Gourdon2003NumericalEO}, \cite{Odlyzko1988Fast}, and \cite{PlattTrudgian2021}.
In contrast, interpolation methods for Dirichlet L-functions have been less explored. However, there are some notable studies on the Riemann zeta function, such as the work by Hiary and Odlyzko \cite{HiaryOdlyzko2012}, who presented an approximation method for the Riemann zeta function using polynomials. Their method, however, requires appropriate scaling of the approximant. Recent advancements by Matiyasevich in 2013 and 2019 introduced an approximation of the Riemann zeta function using a finite Dirichlet series, which has the advantage of not requiring such scaling (see \cite{MATIYASEVICH2020460} and \cite{ matiyasevich2013}). In these works, Matiyasevich, utilizing supercomputers, explored the roles of Gram points and zeros in function approximation and zero discovery through interpolating determinants. 

Our work complements Matiyasevich's work by providing numerical evidence on the role of Gram points and zeros of the Dirichlet L-function \(L(s, \chi)\), corresponding to real Dirichlet characters \(\chi\), as interpolation points. It is important to highlight that naive interpolation, without imposing additional constraints on the approximant, fails to achieve accurate results, necessitating special considerations. After imposing the necessary adjustments, we demonstrate that our method when it leverages \(M\) zeros on the critical line (assuming \textit{the generalized Riemann hypothesis}), it enables the discovery of at least \(M\) additional zeros beyond the \(M\)-th zero. Furthermore, we demonstrate that using \(M\) Gram points \(g_0, \dots, g_M\) enables the approximation of the Dirichlet L-function at \(s = \sigma + it\) for \(|t| \leq g_M\) and slightly beyond, facilitating the discovery of zeros with imaginary parts up to \(g_M\), with errors much smaller than machine epsilon (\(10^{-23}\)).

Our approach does not claim to surpass existing methods; however, on one hand, it offers valuable theoretical insights into the roles of Gram points and zeros of Dirichlet L-functions in function approximation and zero discovery. On the other hand, our approach seeks to inspire further exploration into numerical algorithms for approximating and locating the zeros of Dirichlet \(L\)-functions with very large modulus, about which little to no information is known. In these cases, high numerical instability poses a significant challenge, hence requiring access to supercomputers, which the author lacks. Nevertheless, by utilizing tools from numerical linear algebra and machine learning, specifically feature selection via lasso regression and iterative solvers like the Generalized Minimal Residual and Conjugate Gradient method (see \cite{Saad1986}, \cite{Nifa2017}, and  \cite{hestenes1952} for details), we eliminate the numerical instability problem for real Dirichlet L-functions with small modulus. Hence, since our numerical results demonstrate that it is possible to approximate Dirichlet \(L\)-functions of small modulus without direct function evaluations, this suggests a promising, albeit challenging, path forward, given the computational power of modern supercomputers. For example, by using Gram points as interpolation nodes, points that can be computed without the need for explicit evaluation of the Dirichlet \(L\)-function, there is some hope of discovering some information about Dirichlet L-functions with very large modulus.

\subsection{A Survey of the Riemann Zeta Function and Dirichlet L-functions: Their Zeros and Gram Points}
\subsubsection{The Riemann Zeta Function}

The Riemann zeta function \cite{Ivic1985}, denoted by $\zeta(s)$, is the analytic continuation of the Dirichlet series
\begin{align}
\zeta(s) = \sum_{n=1}^\infty \frac{1}{n^s}, \quad \text{defined for } \Re(s) > 1.
\end{align}
This function is characterized by the functional equation
\begin{align}
\xi(s) = \xi(1-s),
\end{align}
where
\begin{align}
\xi(s) = \frac{1}{2} \pi^{-s/2} s(s-1) \Gamma(s/2) \zeta(s).
\end{align}
From this functional equation, we observe that $\zeta(s)$ has zeros at the negative even integers $-2, -4, -6, \ldots$, known as the trivial zeros of the Riemann zeta function.

The rest of the zeros, which are called the non-trivial zeros of the Riemann zeta function, have garnered considerable attention due to their profound connection to the distribution of prime numbers. These zeros are located within the critical strip \(0 < \Re(s) <1 \) and exhibit symmetry about both the real axis and the critical line $\Re(s) = \frac{1}{2}$. Central to this area of study is the Riemann Hypothesis \cite{Ivic1985}, one of the greatest unsolved problems in mathematics, which conjectures that all non-trivial zeros of $\zeta(s)$ lie on the critical line $s = \frac{1}{2} + it$, where $t \in \mathbb{R}$. Consequently, the $m$-th non-trivial zeros can be expressed as $\rho_{\pm m} = \frac{1}{2} \pm i \gamma_{m}$ with $\gamma_{m} \in \mathbb{R^+}$.

The Riemann zeta function admits a polar decomposition on the critical line:
\begin{align}
\zeta\left(\frac{1}{2} + it\right) = e^{-i\theta(t)} Z(t),
\end{align}
where \( \theta(t) \) is the odd analytic Riemann-Siegel theta function defined as 
\begin{align}
\theta(t) &= \arg(\log \Gamma (1/4 +it/2)) - \frac{\log(\pi) t}{2} \nonumber \\
&= \Im (\log \Gamma (1/4 +it/2)) - \frac{\log(\pi) t}{2},
\end{align}
and \( Z(t) \) is the Hardy-Z function defined by 
\begin{align}
    Z(t) = e^{i \theta(t)}\zeta\left(\frac{1}{2} + it\right),
\end{align} 
both of which are real for real $t$ (see \cite{Edwards1974} and \cite{Gabcke1979}).

\begin{figure}[ht]
 \renewcommand{\thefigure}{1.2.1}
  \centering
  \includegraphics[width=0.8\textwidth]{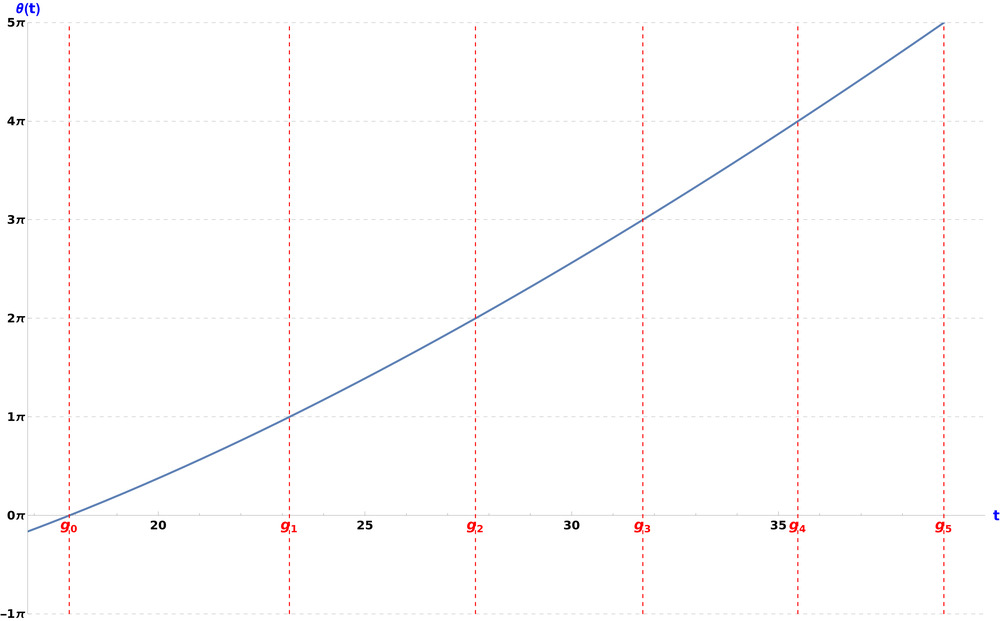}
  \caption{The Riemann-Siegel theta function $\theta(t)$}
  \label{fig:thetaz1}
\end{figure}

A point that satisfies $\theta(t) = m \pi$ for $m \in \mathbb{Z}_{\ge 0}$ is called the $m$-th Gram point and will be denoted $g_m$ (see Figure $1.2.1$). Gram points are particularly interesting since the polar decomposition given in Equation (4) reduces to the following at $t=g_m$ for $m \ge 0$:
\begin{align}
\zeta\left(\frac{1}{2} + ig_m\right) = (-1)^m Z(g_m).
\end{align}

\subsubsection{Dirichlet Characters and the Classification of Real Ones}

Given a natural number \( q \geq 1 \), a Dirichlet character \(\chi\) modulo \( q \) is defined as a group homomorphism \cite{Davenport2000}:
\begin{align}
\chi: (\mathbb{Z}/q\mathbb{Z})^* \rightarrow \mathbb{C}^*.
\end{align}
Since \((\mathbb{Z}/q\mathbb{Z})^*\) is a group of order \(\phi(q)\), where $\phi$ is the Euler's totient function, and \(\chi\) is a group homomorphism, it follows that the image of \(\chi\) is a subset of the unit circle. This character is then extended to all integers by periodicity with period \( q \), and is defined as
\begin{align}
\chi(n) &= 
\begin{cases}
    1 & \text{if } (n,q) = 1, \\ 
    0 & \text{otherwise}.
\end{cases}
\end{align}
The simplest Dirichlet character is the principal Dirichlet character modulo \( q \), defined by
\begin{align}
\chi_0(n) =
\begin{cases}
    1 & \text{if } (n,q) = 1, \\ 
    0 & \text{otherwise}.
\end{cases} 
\end{align}

A Dirichlet character \(\chi\) modulo \(q\) is said to be induced by a character \(\chi^*\) modulo a divisor \(d\) of \(q\) if it satisfies
\begin{align}
\chi(n) = 
\begin{cases}
    \chi^*(n) & \text{if } (n, q) = 1, \\ 
    0 & \text{otherwise}.
\end{cases}
\end{align}
In this context, \(\chi\) can be expressed as \(\chi = \chi^* \chi_0\), where \(\chi_0\) is given by Equation (10). A primitive Dirichlet character is defined as one that is not induced by a character with a smaller modulus. Equivalently, primitive Dirichlet characters modulo \(q\) can be characterized as those with conductor \(q\). The conductor of a character is its smallest quasi-period, where natural number \(d\) is referred to as a quasiperiod of a Dirichlet character \(\chi\) modulo \(q\) if \(\chi(m) = \chi(n)\) whenever \(m \equiv n \pmod{d}\) and \((mn, q) = 1\).

A Dirichlet character \(\chi\) is classified as real if its image is contained within the set of real numbers. Consequently, its image is a subset of \(\{0, \pm 1\}\), as it must be contained in the unit circle. In this paper, we focus on real primitive Dirichlet characters due to a symmetry property (see Section 1.2.3). These characters can be classified as \textit{Kronecker symbols} (see \cite{Davenport2000}).

\subsubsection{Dirichlet L-functions}

The Dirichlet \( L \)-series associated with a Dirichlet character \(\chi\) is given by \cite{Davenport2000} 
\begin{align}
L(s,\chi) = \sum_{n=1}^\infty \frac{\chi(n)}{n^s}, \quad \text{for } \Re(s) > 1.
\end{align}
Due to the complete multiplicativity of \(\chi\), the \( L \)-series admits the Euler product
\begin{align}
L(s,\chi) = \prod_{p \text{ prime}} \left(1 - \frac{\chi(p)}{p^s}\right)^{-1}, \quad \text{for } \Re(s) > 1.
\end{align}

It suffices to study Dirichlet L-series associated to primitive Dirichlet characters because the \( L \)-series for \(\chi\) and \(\chi^*\) relate as follows: \cite{Apostol1976}
\begin{align}
L(s, \chi) =  L(s, \chi^*) \prod_{p \mid q, p \text{ prime}} \left(1 - \frac{\chi^*(p)}{p^s}\right), \quad \Re(s) > 1.
\end{align}

For a primitive Dirichlet character \(\chi\), the \( L \)-series \( L(s,\chi) \) can be analytically continued to an entire function if \(\chi\) is non-principal, or to a meromorphic function with a simple pole at \( s = 1 \) if \(\chi = \chi_0\). The continuation is called the Dirichlet L-function associated to \(\chi\) and also denoted $L(s, \chi)$. The L-function \( L(s,\chi) \) is given by the functional equation \cite{MontgomeryVaughan2006}
\begin{align}
\xi(s) = \epsilon(\chi) \xi(1-s),
\end{align}
where 
\begin{align}
\xi(s) &= \left(\frac{q}{\pi}\right)^{\frac{s+a}{2}} \Gamma\left(\frac{s+a}{2}\right) L(s, \chi), \\
\epsilon(\chi) &= \frac{\tau(\chi)}{i^a \sqrt{q}}.
\end{align}
Here, \(\tau(\chi)\) denotes the Gauss sum
\begin{align}
\tau(\chi) = \sum_{n=1}^q \chi(n) e^{\frac{2 \pi i n}{q}},
\end{align}
and \( a \) is defined as
\begin{align}
a = 
\begin{cases}
    0 & \text{if } \chi(-1) = 1, \\
    1 & \text{if } \chi(-1) = -1.
\end{cases}
\end{align}

Similar to the Riemann zeta function, Dirichlet L-functions have trivial zeros at negative integers, which are either even or odd depending on the sign of \(\chi(-1)\). The remaining non-trivial zeros are located within the critical strip \(0 \le \Re(s) \le 1\) and exhibit symmetry about the critical line \(\Re(s) = \frac{1}{2}\). Additionally, when the associated Dirichlet character is real, the zeros of the L-function are symmetric about the real axis. The Generalized Riemann Hypothesis posits that all zeros of Dirichlet L-functions lie on the critical line \(\Re(s) = \frac{1}{2}\) (see \cite{Davenport2000}). 

For a primitive Dirichlet character \(\chi\), analogous to the Riemann Zeta function, we can express \( L(s,\chi) \) in polar form as 
\begin{align}
L\left(\frac{1}{2} + it,\chi\right) = e^{-i \theta(t, \chi)} Z(t,\chi) \quad (\text{see } [\cite{Rumely1993}, \text{ p. } 426]),
\end{align}
where $\theta(t, \chi)$ and $Z(t,\chi)$ are real for real $t$, and they are defined as follows
\begin{align}
\theta(t, \chi) &= \arg\left(\log \Gamma \left(\frac{1}{4} + \frac{a}{2} + \frac{it}{2}\right)\right) - \frac{\log(\pi) t}{2} + \frac{\log(q) t}{2} + \frac{i}{2} \log(\epsilon(\chi)) \nonumber \\
&= \Im \left(\log \Gamma \left(\frac{1}{4} + \frac{a}{2} + \frac{it}{2}\right)\right) - \frac{\log(\pi) t}{2} + \frac{\log(q) t}{2} + \frac{i}{2} \log(\epsilon(\chi)),
\end{align}
and
\begin{align}
Z(t, \chi) = e^{i \theta(t, \chi)} L\left(\frac{1}{2} + it, \chi\right).
\end{align}

\textbf{Definition 1.2.3.1:} The classical Gram points, denoted \(g_m\) for \(m \in \mathbb{Z}_{\ge 0}\), are the real solutions to \(\theta(t) = m \pi\) for \(m \in \mathbb{Z}_{\ge 0}\). We define the \textit{generalized Gram points} associated with a primitive Dirichlet character \(\chi\) as the real solution(s) to \(\theta(t, \chi) = m \pi\) for \(m \in \mathbb{Z}_{\ge 0}\). We denote these points by \(g_m\) for \(m \in \mathbb{Z}_{\ge 0}\) with the associated character \(\chi\) being understood from the context.

\section{Interpolating Dirichlet L-functions}

\subsection{The Numerical Method}

Consider a Dirichlet series with real coefficients \(\{l_n\}\):
\begin{align}
    L(s) = \sum_{n=1}^\infty \frac{l_n}{n^s},
\end{align}
which can be meromorphically continued to the complex plane \(\mathbb{C}\) via a functional equation. Assume \(L(s)\) has a polar decomposition, similar to Equations (4) and (22), and satisfies a Riemann Hypothesis, i.e., its zeros lie on the line \(\Re(s) = \frac{1}{2}\). Suppose we are interested in approximating and locating the zeros of \(L(s)\).

When attempting to interpolate \(L(s)\) using a finite Dirichlet series
\begin{align}
    F_{a_n,~ N}  = \sum_{n = 1}^N \frac{a_n}{n^s},
\end{align}
three questions arise:

\begin{enumerate}[label=\arabic*.]
    \item What distribution should the interpolating points follow to ensure an accurate approximation of the function?
    \item Are the values of the function \(L(s)\) accessible at the chosen interpolating points?
    \item What numerical challenges might arise when solving the system of equations required for interpolation?
\end{enumerate}

Our numerical discoveries indicate that interpolating \( L \) at its zeros or similar distributed points tends to yield accurate approximations when the approximant is a finite Dirichlet series \( F_{a_n,~ N}  \) with additional constraints. Specifically, it is crucial to require \( \{a_n\} \) to consist of real numbers, analogous to \( \{l_n\} \), in order to get reasonable coefficients and approximations. For the Riemann zeta function, the constraints include the latter constraint and the requirement that $a_1= l_1 =1$ (see \cite{MATIYASEVICH2020460, matiyasevich2013}). However, for a Dirichlet L-function corresponding to a Dirichlet character $\chi$, additional constraints are necessary. These include:
\begin{itemize}[label=--]
    \item The coefficients \(a_n\) should be equal to \(l_n = 0\) when \((n, q) > 1\). This condition is suggested by a feature selection experiment conducted by the author.
    \item We should enforce at least \(a_1 = l_1 = \chi(1) = 1\). For example, when interpolating at Gram points, it is necessary to assume \(a_j = l_j = \chi(j)\) for \(j \leq 2\).
\end{itemize}

In numerical approximation, interpolation problems often face ill-conditioning, where small changes in inputs lead to disproportionately large changes in outputs (see \cite{Trefethen2021} and \cite{DeBoor2001}). In our problem, ill-conditioning arises because the interpolating nodes are noisy in the sense that we do not know the exact zeros or Gram points but only approximate values of them. Mitigating this issue typically involves high-precision computation on supercomputers or employing suitable iterative methods. The iterative solver we employ, developed by Saad and Schultz in 1986 (see \cite{Saad1986}), is well-suited for solving interpolation systems. Our first-attempt numerical method is structured as follows:

\begin{enumerate}
    \item Choose the number of coefficients \( N = 2M + k \) to use in the approximant $F_{a_n,~ N} $, where \( k \) determines the number of initial coefficients \( a_n \) tha are forced to match \( l_n \). For example, when \( k = 2 \), the approximant $F_{a_n,~ N} $ should have \( a_1 = l_1  \) and \( a_2 = l_2 \).
    \item Select interpolation points \( x_m \) for \( m = 1, 2, \ldots, M \). We choose these points as initial zeros of \( L \), Gram points of \( L \), or similar distributed points.
    \item Assuming the coefficients $a_n$ are real, formulate the interpolation problem: 
   \begin{align}
    & F_{a_n,~ N} (1/2 + i x_m) = L(1/2 + i x_m)  \quad \text{for} \quad m = 1, \ldots, M, \label{eq:interpolation1} \\
    & \iff \ \Re(F_{a_n,~ N} (1/2 + i x_m)) = \Re(L(1/2 + i x_m))  \label{eq:interpolation2} \\
    & \quad \text{and} \quad \Im(F_{a_n,~ N} (1/2 + i x_m)) = \Im(L(1/2 + i x_m)) \quad \text{for} \quad m = 1, \ldots, M. \label{eq:interpolation3} \\
    & \iff \sum_{n = k+1}^N a_n \Re(n^{-1/2- i x_m}) = \Re(L(1/2 + i x_m)) - \sum_{n=1}^k a_n \Re(n^{-1/2- i x_m})  \\
    & \quad \text{and} \quad \sum_{n = k+1}^N a_n \Im(n^{-1/2- i x_m}) = \Im(L(1/2 + i x_m)) - \sum_{n=1}^k a_n \Im(n^{-1/2- i x_m})  \\
    & \quad \quad \quad \quad \quad \text{for} \quad m = 1, \ldots, M. \nonumber
\end{align}
The reason for solving the interpolation problem defined by Equations (28) and (29), rather than the equivalent problem in Equation (25), is to ensure that the solution to the system is a vector of \textbf{real} coefficients \(a_n\).
    \item Solve the system 
    \begin{align}
        A x = b 
    \end{align} iteratively, where \( A = (\mu_{m,n}) \) is an \( (N-k) \times (N-k) \) matrix, $x$ is an $N-k$ vector of the unknown coefficients $a_n$ for $ k+1\le n \le N$, and $b=(b_m) $ is an $N-k$ vector. From step 3, we have that $\mu_{m,n}$ and $b_m$ are defined as follows:
    \begin{align}
    \mu_{m,n} =
    \begin{cases}
    \Re(n^{-1/2 - i x_m}) & \text{if } 1 \leq m \leq M, \\
    \Im(n^{-1/2 - i x_m}) & \text{if } M <m \leq 2M,
    \end{cases} \label{eq:system_matrix}
    \end{align}
    and 
    \begin{align}
    b_m =
    \begin{cases}
    \Re(L(1/2 + i x_m)) - \sum_{n=1}^k a_n \Re(n^{-1/2 - i x_m}) & \text{if } 1 \leq m \leq M, \\
    \Im(L(1/2 + i x_m)) - \sum_{n=1}^k a_n \Im(n^{-1/2 - i x_m}) & \text{if } M <m \leq 2M.
    \end{cases} \label{eq:system_rhs}
    \end{align}
    \item Step $4$ determines \( a_n \) for \(  k+1 \leq n \leq N \), and we set \( a_n = l_n \) for \( n \leq k \). The resulting approximant is:
    \begin{align}
    F_{a_n,~ N} (s) = \sum_{n = 1}^N \frac{a_n}{n^s}. \label{eq:approximant}
    \end{align}
\end{enumerate}

In Step 4, the author constructed and solved the linear system presented in Equation (30) using Mathematica. For this, the Generalized Minimal Residual Method (GMRES) was employed, as described in \cite{Saad1986}. Below is the Mathematica code used to solve the system with the GMRES method:

\begin{lstlisting}{mathematica}
x = LinearSolve[
  A, b, 
  Method -> {"Krylov", Method -> "GMRES", 
    Tolerance -> 10^-20, MaxIterations -> 1000}
];
\end{lstlisting}

In this code snippet, \texttt{Tolerance} specifies the tolerance for the relative residual norm \( \frac{\|b - A x\|}{\|b\|} \leq \text{tol} \), with a default of \( 1 \times 10^{-6} \), and \texttt{MaxIterations} limits the maximum number of iterations. 

\subsection{Utilizing Zeros of Dirichlet L-functions to Discover New Zeros}

Numerical experiments have shown that leveraging \( M > 0 \) initial zeros of a Dirichlet \( L \)-function $L(s, \chi)$, associated with a real Dirichlet character \(\chi\) of modulus \( q \) and a fundamental discriminant \( d \),
\begin{align}
    L(s, \chi) = \sum_{n=1}^\infty \frac{\chi(n)}{n^s},
\end{align}
enables the discovery of at least \( M \) additional zeros beyond the \( M \)-th known zero. However, as previously mentioned, solving the system in Equation (30) without improsing additional constraints on the approximant \( F_{a_n,~ N}  \) is insufficient to guarantee that \( F_{a_n,~ N}  \) serves as an effective approximant. Specifically, when \( k=0 \) and \( x_m = \rho_m \) for \( m=1,2,\ldots,M \) (where \( \rho_m \) denotes the \( m \)-th zero of the corresponding Dirichlet \( L \)-function), \( b_m = 0 \) for \( m=1,2,\ldots,2M \) (see Equation (32)). This results in a homogeneous system of equations (Equation (30)), leading to a trivial solution. To avoid such a trivial result, we consider \( k=1 \) (thus \( N = 2M +1 \)), by imposing \( a_1 = l_1 = \chi_q(1) = 1 \). Consequently, Equation (32) simplifies to:
\begin{align}
    b_m =
    \begin{cases}
    -1, & \text{if } 1 \leq m \leq M, \\
    0, & \text{if } M < m \leq 2M.
    \end{cases} \label{eq:system_rhs2}
\end{align}
\subsubsection{Case of \( q = 4 \) and \( d = -4 \)}
As discussed, solving the interpolation problem in Equation (30), without more constraints on the approximant \( F_{a_n,~ N} \), gives ineffective results. For instance, with \( q=4 \), \( d=-4 \), and \( M=200 \), the coefficients obtained from the system of equations render \( F_{a_n,~ 401} \) ineffective. For example, the value of \( F_{a_n,~ 401} \) at the first interpolating point \( \rho_1 \), which is supposed to be zero, is approximately \( 2.0643939 \dots + 0.9867209 \dots i \), which is absurd. Additionally, the coefficients exhibit problematic behavior, as illustrated in Figure 2.2.1.1.

\begin{figure}[ht]
 \renewcommand{\thefigure}{2.2.1.1}
  \centering
  \includegraphics[width=0.7\textwidth]{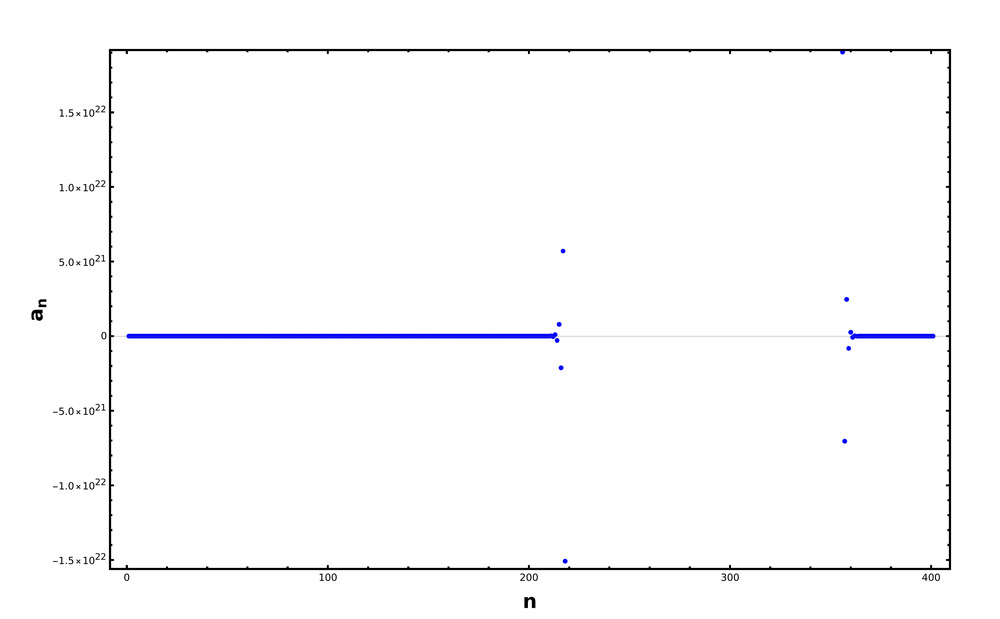}
  \caption{Problematic behavior of the coefficients \( a_n \) of the approximant \( F_{a_n,~ 401}  \).}
  \label{fig:thetaz2}
\end{figure}

A partial explanation for this failure is that \( F_{a_n,~ N}  \) does not sufficiently share the fundamental properties of the Dirichlet \( L \)-function \( L(s,\chi_4) \). One approach to address this issue is to ensure that the coefficients of the approximant exhibit behavior similar to that of \( \chi_4(n) \) in \( L_{\chi_4} \). To address this issue, we conducted a machine learning experiment by treating the coefficients of the approximant as features in the language of machine learning and by applying lasso regression to identify the most important features (coefficients). The coefficients that vanished last during the regression process were deemed the most important for determining the approximant (see \cite{fonti2017} for details on how feature selection is applied using lasso regression). Based on the conducted experiment, we concluded that the coefficients \(a_n\) for \((n,4) > 1\) are the least significant features. Consequently, the condition \(a_n = \chi(n) = 0\) was imposed for \((n,4) > 1 \iff (n,2) > 1\). Accordingly, we adjusted the system of equations (30) by setting \( A = (\mu_{m, 2n+1})_{1 \leq m,n \leq 2M} \) and defining \( b \) by Equation (35). The approximant \( F_{a_n,~ 401} \) is now given by Figure 2.2.1.2.
\begin{figure}[ht]
 \renewcommand{\thefigure}{2.2.1.2}
  \centering
  \includegraphics[width=0.7\textwidth]{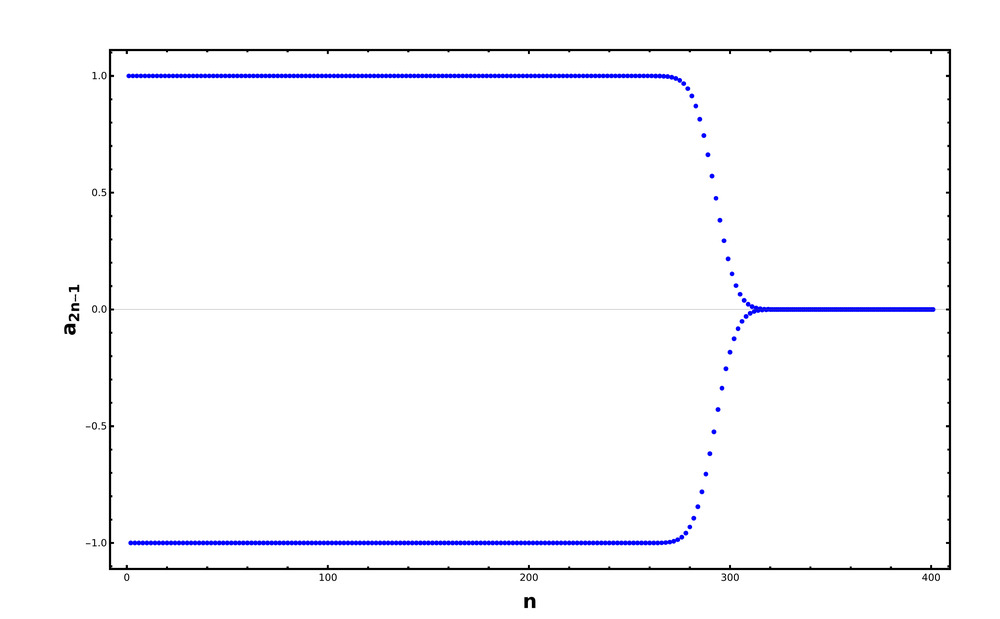}
  \caption{Coefficients \( a_{2n+1} \) of the approximant \( F_{a_n,~ 401}  \) when \( M=200 \).}
  \label{fig:thetaz3}
\end{figure}

\vspace{5 cm}
Empirical evidence demonstrates the effectiveness of this modified approximant. The approximant \( F_{a_n,~ 401}  \) provides surprisingly accurate approximations for more than 200 zeros of \( L(s,{\chi_{4}}) \) beyond the 200th zero, as illustrated in Equations (36) to (44).

\begin{align}
0 &= F_{a_n,~ 401}  \left( \rho_{201} + 9.75349 \dots \times 10^{-133} - 4.18828 \dots \times 10^{-132} i \right) \\
0 &= F_{a_n,~ 401}  \left( \rho_{202} + 7.32343 \dots \times 10^{-131} - 6.78124 \dots \times 10^{-131} i \right) \\
& \dots \nonumber \\ 
0 &= F_{a_n,~ 401}  \left( \rho_{250} + 3.27918 \dots \times 10^{-81} - 9.2487 \dots \times 10^{-81} i \right) \\
& \dots \nonumber \\
0 &= F_{a_n,~ 401}  \left( \rho_{300} - 1.04 \dots \times 10^{-54} - 2.49132 \dots \times 10^{-54} i \right) \\
& \dots \nonumber \\ 
0 &= F_{a_n,~ 401}  \left( \rho_{350} - 8.23308 \dots \times 10^{-37} - 2.92157 \dots \times 10^{-37} i \right) \\
& \dots \nonumber \\
0 &= F_{a_n,~ 401}  \left( \rho_{400} + 1.76021 \dots \times 10^{-25} + 2.18751 \dots \times 10^{-25} i \right) \\
& \dots \nonumber \\ 
0 &= F_{a_n,~ 401}  \left( \rho_{450} - 2.29465 \dots \times 10^{-17} + 5.12563 \dots \times 10^{-18} i \right) \\
& \dots \nonumber \\ 
0 &= F_{a_n,~ 401}  \left( \rho_{500} - 8.8619 \dots \times 10^{-11} - 1.99089 \dots \times 10^{-10} i \right) \\
& \dots \nonumber \\ 
0 &= F_{a_n,~ 401}  \left( \rho_{510} - 3.84595 \dots \times 10^{-10} + 3.14449 \dots \times 10^{-10} i \right) \\
& \dots  \nonumber
\end{align}

We may use the 200 discovered zeros as interpolating points to further locate additional zeros beyond the 400th zero. Equations (45) to (51) reveal the performance of the approximant \( F_{a_n,~ 801}  \), utilizing \( M = 400 \) zeros:

\begin{align}
0 &= F_{a_n,~ 801}  \left( \rho_{401} + 1.20485 \dots \times 10^{-140} - 7.73352 \dots \times 10^{-141} \, \mathrm{i} \right) \\
& \dots \nonumber \\ 
0 &= F_{a_n,~ 801}  \left( \rho_{500} + 4.60093 \dots \times 10^{-85} - 1.52913 \dots \times 10^{-82} \, \mathrm{i} \right) \\
& \dots \nonumber \\ 
0 &= F_{a_n,~ 801}  \left( \rho_{600} + 4.58311 \dots \times 10^{-83} - 5.25731 \dots \times 10^{-83} \, \mathrm{i} \right) \\
& \dots \nonumber 
\end{align}
\begin{align} 
0 &= F_{a_n,~ 801}  \left( \rho_{700} + 3.8944 \dots \times 10^{-81} + 1.76879 \dots \times 10^{-81} \, \mathrm{i} \right) \\
& \dots \nonumber \\
0 &= F_{a_n,~ 801}  \left( \rho_{800} + 2.22741 \dots \times 10^{-54} + 2.70752 \dots \times 10^{-53} \, \mathrm{i} \right) \\
& \dots \nonumber \\
0 &= F_{a_n,~ 801}  \left( \rho_{900} + 7.80376 \dots \times 10^{-37} + 2.97914 \dots \times 10^{-35} \, \mathrm{i} \right) \\
& \dots \nonumber \\ 
0 &= F_{a_n,~ 801}  \left( \rho_{1000} + 2.28595 \dots \times 10^{-21} + 1.09729 \dots \times 10^{-21} \, \mathrm{i} \right) \\
& \dots \nonumber
\end{align}
\subsubsection{Case of \( q = 3 \) and \( d = -3 \)}
Fix \( M \) and \( k \) as positive integers. Similar to the approach in Section 2.2.1, we must ensure that \( a_n = 0 \) for any \( n \) where \((n, 3) > 1\). To achieve this, we modify the system of equations (30) by setting \( A = (\mu_{m, n})_{1 \leq m \leq 2M, n \in I \setminus J} \) and defining \( b \) as in Equation (35). Here, \( I \) is a list of length \( N = 2M + k \) that contains the first \( 2M + k \) natural numbers \( n \) with \((n, 3) = 1\) (i.e., non-multiples of 3). The set \( J \) includes the first \( k \) non-multiples of 3, representing the indices of the coefficients \( a_n \) in \( F_{a_n, N} \) that are fixed to match the corresponding coefficients of \( L(s, \chi_3) \), meaning \( a_n = \chi_3(n) \) for \( n \in J \).

For instance, by setting \( M = 200 \) and \( k = 1 \), the resulting approximant \( F_{a_n,~ 401} \) is illustrated in Figure 2.2.2.2. Similarly, the approximant \( F_{a_n,~ 401}  \) provides surprisingly accurate approximations for more than 200 zeros of \( L(s,{\chi_{3}}) \) beyond the 200th zero, as illustrated in Equations (52) to (58). For example:

\begin{figure}[ht]
 \renewcommand{\thefigure}{2.2.2.2}
  \centering
  \includegraphics[width=0.7\textwidth]{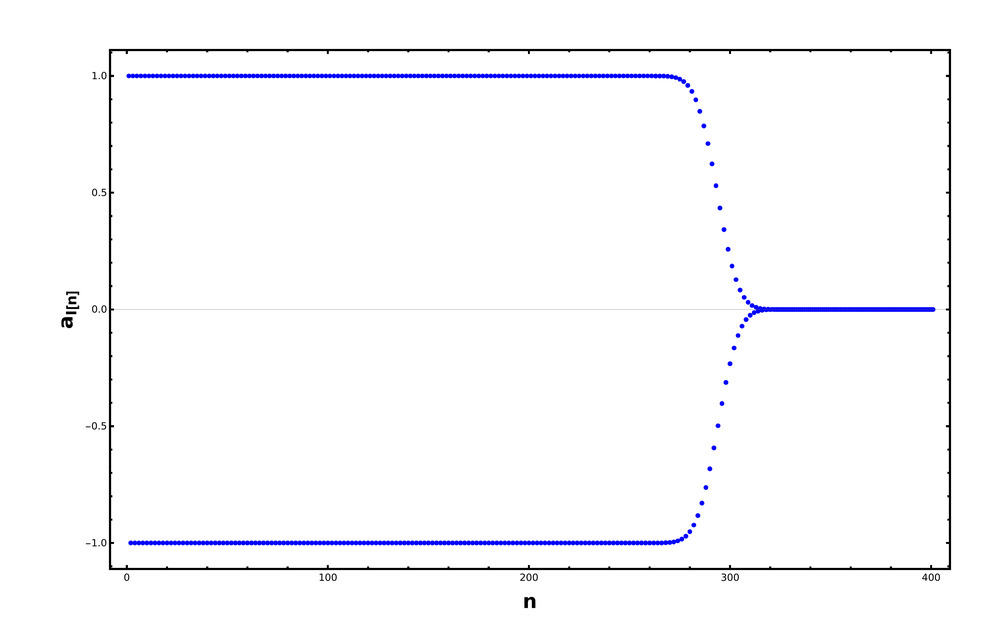}
  \caption{Coefficients \( a_{I[n]} \) of the approximant \( F_{a_n,~ 401}  \) when \( M=200 \), where \( I[n] \) is the \( n \)-th element of the list $I$}
  \label{fig:thetaz4}
\end{figure}

\begin{align} 
0 &= F_{a_n,~ 401}  \left( \rho_{201} + 3.096 \dots \times 10^{-126} -1.89975 \dots \times 10^{-126} i \right) \\
& \dots \nonumber \\ 
0 &= F_{a_n,~ 401}  \left( \rho_{250} -2.12089 \dots \times 10^{-82} +1.19414 \dots \times 10^{-82} i \right) \\
& \dots \nonumber \\
0 &= F_{a_n,~ 401}  \left( \rho_{300} + 2.17051 \dots \times 10^{-51} -4.35149 \dots \times 10^{-51} i \right) \\
& \dots \nonumber \\  
0 &= F_{a_n,~ 401}  \left( \rho_{350} -3.91016 \dots \times 10^{-38} +1.03771 \dots \times 10^{-35} i \right) \\
& \dots \nonumber \\ 
0 &= F_{a_n,~ 401}  \left( \rho_{400} -1.12539 \dots \times 10^{-23} -1.14114 \dots \times 10^{-23} i \right) \\
& \dots \nonumber \\  
0 &= F_{a_n,~ 401}  \left( \rho_{450} -1.40184 \dots \times 10^{-19} -1.91713 \dots \times 10^{-19} i \right) \\
& \dots \nonumber \\  
0 &= F_{a_n,~ 401}  \left( \rho_{500} + 4.34172  \dots \times 10^{-12} -2.87374 \dots \times 10^{-10} i \right) \\
& \dots  \nonumber
\end{align}

As before, we may use the 200 discovered zeros as interpolating points to further locate additional zeros beyond the 400th zero. Equations (59) to (65) reveal the performance of the approximant \( F_{a_n,~ 801}  \), utilizing \( M = 400 \) zeros:

\begin{align}   
0 &= F_{a_n,~ 801}  \left( \rho_{401} + 1.78911 \dots \times 10^{-142} -8.69534 \dots \times 10^{-141} \, \mathrm{i} \right) \\
& \dots \nonumber \\ 
0 &= F_{a_n,~ 801}  \left( \rho_{500} + 9.96048 \dots \times 10^{-96} -1.36208 \dots \times 10^{-96} \, \mathrm{i} \right) \\
& \dots \nonumber \\ 
0 &= F_{a_n,~ 801}  \left( \rho_{600} + 3.22705 \dots \times 10^{-96} -8.83612 \dots \times 10^{-96} \, \mathrm{i} \right) \\
& \dots \nonumber \\ 
0 &= F_{a_n,~ 801}  \left( \rho_{700} + 3.11343 \dots \times 10^{-77} -2.88703 \dots \times 10^{-77} \, \mathrm{i} \right) \\
& \dots \nonumber \\ 
0 &= F_{a_n,~ 801}  \left( \rho_{800} + 2.93314 \dots \times 10^{-51} -3.33836 \dots \times 10^{-51} \, \mathrm{i} \right) \\
& \dots \nonumber 
\end{align}
\begin{align}
0 &= F_{a_n,~ 801}  \left( \rho_{900} + 5.77277 \dots \times 10^{-35} -2.16093 \dots \times 10^{-34} \, \mathrm{i} \right) \\
& \dots \nonumber \\ 
0 &= F_{a_n,~ 801}  \left( \rho_{1000} + 7.78242 \dots \times 10^{-22} + 3.74379 \dots \times 10^{-21} \, \mathrm{i} \right) \\
& \dots \nonumber
\end{align}

\subsection{A Yet Less Expensive and Comparitively as Efficient Approach for Approximating and Locating Zeros of Dirichlet L-functions}

Numerical experimentation showed that one can impose fewer constraints imposed on \( F_{a_n,~ N}  \) and yet still maintain a surprisingly good approximant of a Dirichlet \( L \)-function $L(s, \chi)$ with $\chi$ a real primitive Dirichlet character. Specifically, we may drop half of the equations provided by Equations (28) and (29). The advantage of this approach, provided it yields satisfactory approximation results, is that one may replace the use of zeros as interpolating points with other points of similar density and spacing, which are easier to compute and still do not require the evaluation of $L(s, \chi)$ at the interpolating points. 

One trivial inefficient approach to selecting new interpolating points of similar spacing might be to use \( \rho_m + 1 \) for \( m = 1, 2, \ldots, M \). However, this choice is counterproductive if the goal is to minimize function evaluations of \( L(s, \chi) \) and reduce computational costs. Instead, the aim is to identify alternative interpolating points that are computationally less expensive, ideally where \( \Re\left(L\left(\frac{1}{2} + i t\right)\right) \) or \( \Im\left(L\left(\frac{1}{2} + i t\right)\right) \) vanishes.

Since \( Z(t, \chi) \) is real for real \( t \), we have
\begin{align}
    \Im\left( e^{i \theta(t, \chi)} L\left(\frac{1}{2} + it, \chi\right)\right) = 0 \quad \text{for real } t.
\end{align}
Thus, one can impose the following constraint for some real interpolating points \( \{x_m\}_{m=1}^M \), by requiring 
\begin{align}
    \Im\left( e^{i \theta(x_m, \chi)} L\left(\frac{1}{2} + ix_m, \chi\right)\right) = 0 \quad \text{for } m = 1, \dots, M.
\end{align}
However, we may further simplify the last equation by taking \( x_m = g_m \) for \( m = 1, \dots, M \), where \( g_m \) are the generalized Gram points of \( L(s, \chi) \). Since \( e^{i \theta(g_m, \chi)} = (-1)^m \), we can simplify the system of equations to
\begin{align}
    \Im\left( L\left(\frac{1}{2} + ig_m, \chi\right)\right) = 0 \quad \text{for } m = 1, \dots, M.
\end{align}

Now choose positive natural numbers \( M \) and \( k \), and let \( I \) be the list of \( M + k \) natural numbers \( n \) such that \( (n, q) > 1 \), where \( q \) is the modulus of the character. Let \( J \) be the sublist of \( I \) containing the first \( k \) natural numbers of \(I\). The elements of \( J \) are the indices of the coefficients \( a_n \) of \( F_{a_n,~ N}  \) that must match the coefficients of \( L(s,\chi_q) \) with the same indices, i.e., \( a_n = \chi(n)\) for $n \in J$. The linear system \( A x = b \) in Equation (30) now becomes of size \( (N - k) \times (N - k) \), where \( N \) is now taken to be \( M + k \) instead of \( 2M + k \). Recall that \( M \) is the number of Gram points used for the interpolation task. 

However, just as in Section 2.2, the approximant does not yet ``know enough'' about \( L(s,\chi_q) \), so we force the coefficients \( a_n \) of \( F_{a_n,~ N}  \) to satisfy \( a_n = 0 \) for \( (n, q) > 1 \). In other words, the system becomes \( A = (\mu_{m,n})_{1 \le m \le M, n \in I \setminus J} \) and \( b = (b_m)_{1 \le m \le M} \), where
\begin{align}
    \mu_{m,n} = \Im\left(n^{-1/2 - i g_m}\right) \quad \text{for } 1 \le m \leq M, \, n \in I \setminus J
\end{align}
and 
\begin{align}
    b_m = - \sum_{n \in J} \chi(n) \Im\left(n^{-1/2 - i g_m}\right) \quad \text{for } 1 \le m \leq M.
\end{align}
The resulting approximant is:
\begin{align}
    F_{a_n,~ N} (s) = \sum_{n \in I} \frac{a_n}{n^s}. \label{eq:approximant2}
\end{align}

Forcing only one coefficient of \( F_{a_n,~ N}  \) to match the first coefficient of \( L(s,{\chi_q}) \) is no longer sufficient. This is because, for \( k = 1 \), the right-hand side of Equation (70) becomes the zero quantity $\chi_q(1) \Im\left(1^{-1/2 - i g_m}\right)$, since the list \( I \) always contains 1 as its first element, resulting in a trivial solution. To address this issue, we adjust by considering \( k = 2 \) instead.

\subsubsection{Case of \( q = 4 \) and \( d = -4 \)}
When \( q = 4 \) and \( d = -4 \), the Gram points \( \{ g_m \} \) are illustrated in Figure $2.3.1.1$.  Given that \( q = 4 \) and \( k = 2 \), we set the coefficients of even indices to zero, and we set \( a_1 = \chi_4(1) = 1 \) and \( a_3 = \chi_4(3) = -1 \). Utilizing \( M = 500 \) generalized Gram points (and thus \( N = 502 \)) to find the approximant \( F_{a_n,~ 502}  \), the coefficients of \( F_{a_n,~ 502}  \) are specified by Figure $2.3.1.2.$

\begin{figure}[ht]
 \renewcommand{\thefigure}{2.3.1.1}
  \centering
  \includegraphics[width=0.7\textwidth]{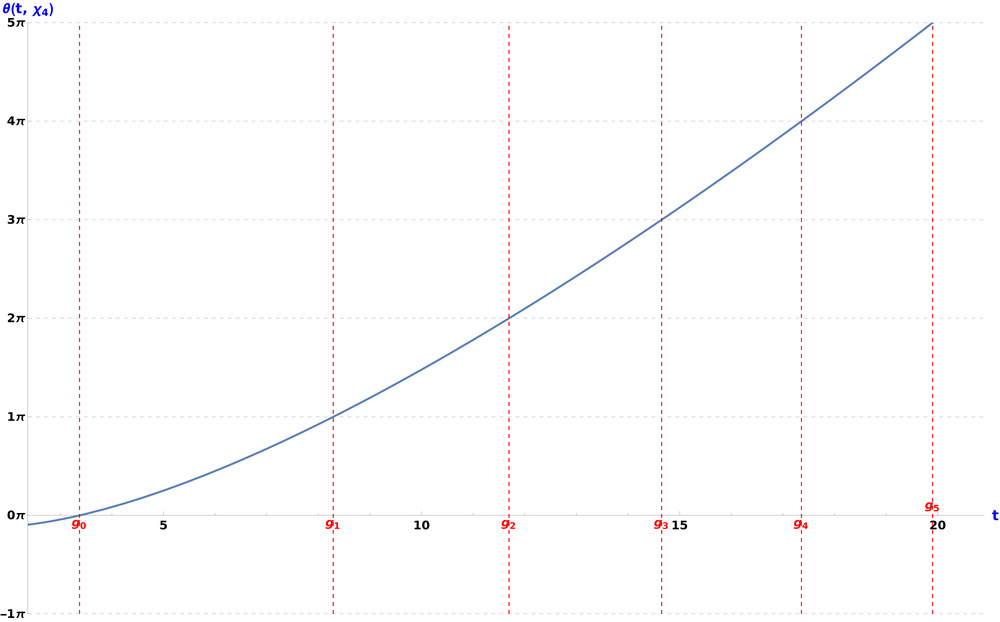}
  \caption{The Riemann-Siegel theta function $\theta(t, \chi_4)$}
  \label{fig:thetaz5}
\end{figure} 

\begin{figure}[ht]
 \renewcommand{\thefigure}{2.3.1.2}
  \centering
  \includegraphics[width=0.7\textwidth]{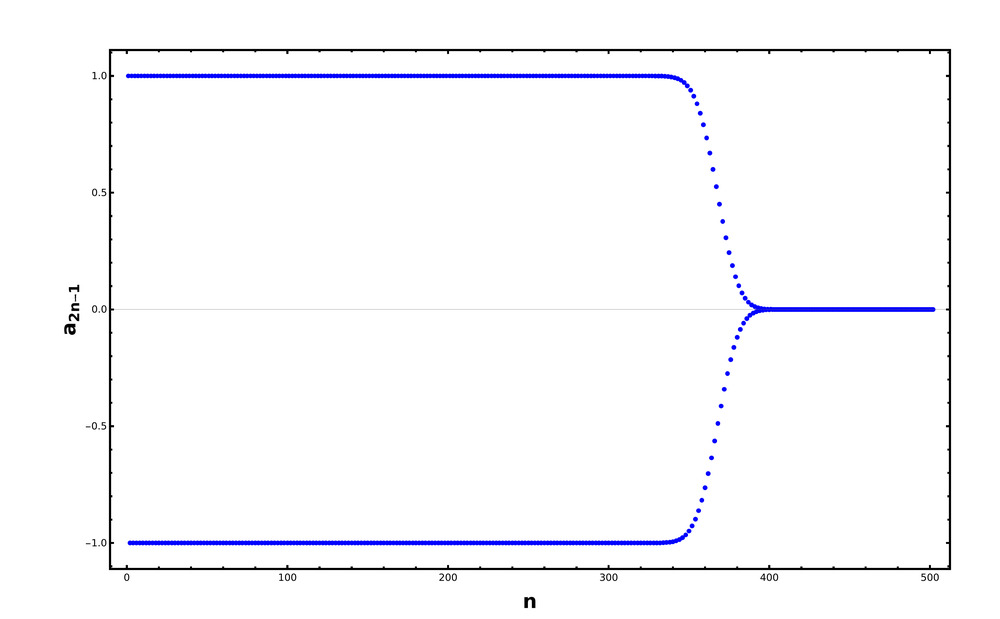}
  \caption{Coefficients \( a_{2n-1} \) of the approximant \( F_{a_n,~ 502}  \) for \( M=500 \).}
  \label{fig:thetaz6}
\end{figure} 
 
\begin{table}
\centering
\caption{Approximation of $ L(s,{\chi_4})$ by $F_{a_n,~ 502}$.}

\begin{tabular}{cc}
    \begin{subtable}[t]{0.45\textwidth}
    \centering
    \begin{tabular}{c c}
    \toprule
    $s$ & \multicolumn{1}{c}{$|L(s,{\chi_4})- F_{a_n,~ 502} (s)|$} \\
    \midrule
    $-4 + 100i$ & $5.7969 \dots \times 10^{-87}$ \\
    $-3 + 300i$ & $1.73264 \dots \times 10^{-71}$ \\
    $-2 + 500i$ & $2.64970 \dots \times 10^{-46}$ \\
    $-1 + 640i$ & $5.86596 \dots \times 10^{-33}$ \\
    $-1/2 + 100i$ & $9.62948 \dots \times 10^{-95}$ \\
    $-1/2 + 300i$ & $3.42636 \dots \times 10^{-77}$ \\
    $-1/2 + 500i$ & $4.66234 \dots \times 10^{-50}$ \\
    $-1/2 + 640i$ & $2.90574 \dots \times 10^{-34}$ \\
    $100i$ & $9.02855 \dots \times 10^{-96}$ \\
    $300i$ & $2.47881 \dots \times 10^{-78}$ \\
    \bottomrule
    \end{tabular}
    \end{subtable}
    &
    \begin{subtable}[t]{0.45\textwidth}
    \centering
    \begin{tabular}{c c}
    \toprule
    $s$ & \multicolumn{1}{c}{$|L(s,{\chi_4})- F_{a_n,~ 502} (s)|$} \\
    \midrule
    $500i$ & $2.61296 \dots \times 10^{-51}$ \\
    $640i$ & $1.43945 \dots \times 10^{-35}$ \\
    $1/2 + 100i$ & $9.70096 \dots \times 10^{-97}$ \\
    $1/2 + 300i$ & $1.79354 \dots \times 10^{-79}$ \\
    $1/2 + 500i$ & $1.46451 \dots \times 10^{-52}$ \\
    $1/2 + 640i$ & $7.13111 \dots \times 10^{-37}$ \\
    $1 + 100i$ & $1.91335 \dots \times 10^{-97}$ \\
    $2 + 300i$ & $6.79948 \dots \times 10^{-83}$ \\
    $3 + 500i$ & $8.10845 \dots \times 10^{-59}$ \\
    $4 + 640i$ & $5.22922 \dots \times 10^{-46}$ \\
    \bottomrule
    \end{tabular}
    \end{subtable}
\end{tabular}
\end{table}

The approximant \( F_{a_n,~ 502} \) provides a highly accurate approximation of the zeros of \( L(s,{\chi_4}) \), particularly within the range of imaginary parts from \( g_0  \) to \( g_{499} \), and even slightly beyond \( g_{499} \), as demonstrated in Equation (72) to (79). Remarkably, \( F_{a_n,~ 502} \) was constructed without prior knowledge of any of the zeros of \( L(s, \chi_4) \), yet it detects zeros on the critical line with imaginary parts between \( \Im(\rho_1) = 6.020948 \dots \) and \( \Im(\rho_{450}) = 628.824833 \dots \), and slightly beyond that.

\begin{align} 
0 &= F_{a_n,~ 502}  \left( \rho_{1} -4.69865 \dots \times 10^{-91} -4.43482 \dots \times 10^{-91} i \right) \\
0 &= F_{a_n,~ 502}  \left( \rho_{2} +  7.00017 \dots \times 10^{-91} +2.21661 \dots \times 10^{-91}  i \right) \\
& \dots \nonumber \\ 
0 &= F_{a_n,~ 502}  \left( \rho_{100} -5.59356 \dots \times 10^{-98} +8.54232 \dots \times 10^{-98}  i \right) \\
& \dots \nonumber \\ 
0 &= F_{a_n,~ 502}  \left( \rho_{200} +1.35254 \dots \times 10^{-81} +2.71801 \dots \times 10^{-80} i  \right) \\
& \dots \nonumber \\ 
0 &= F_{a_n,~ 502}  \left( \rho_{300} +1.52975 \dots \times 10^{-62} +4.42073 \dots \times 10^{-62} i  \right) \\
& \dots \nonumber \\ 
0 &= F_{a_n,~ 502}  \left( \rho_{400} -4.74778 \dots \times 10^{-47} +4.2118 \dots \times 10^{-47} i  \right) \\
& \dots \nonumber \\ 
0 &= F_{a_n,~ 502}  \left( \rho_{460} + 7.52525 \dots \times 10^{-39} +7.32494 \dots \times 10^{-39} i  \right) \\
& \dots \nonumber \\ 
0 &= F_{a_n,~ 502}  \left( \rho_{500} + 5.16403 \dots \times 10^{-33} +1.03997 \dots \times 10^{-33} i \right) \\
& \dots \nonumber
\end{align}

Table 2 provides a sample of the behavior of the approximant \( F_{a_n,~ 702} \), using the interpolating Gram points \( g_0 = 3.3697 \ldots \), \(\ldots\), and \( g_{699} = 831.7394 \ldots \), which results in an improved approximation with reduced error terms.
\begin{table}[ht]
\centering
\caption{Approximation of $L(s, \chi_4)$ by $F_{a_n,~ 702}$.}

\begin{tabular}{cc}
    \begin{subtable}[t]{0.45\textwidth}
    \centering
    \begin{tabular}{c c}
    \toprule
    $s$ & \multicolumn{1}{c}{$|L(s, \chi_4) - F_{a_n,~ 702}(s)|$} \\
    \midrule
    $-4 + 100i$ & $1.178 \cdots \times 10^{-131}$ \\
    $-3 + 300i$ & $3.695 \cdots \times 10^{-127}$ \\
    $-2 + 500i$ & $5.508 \cdots \times 10^{-96}$ \\
    $-1 + 600i$ & $2.853 \cdots \times 10^{-84}$ \\
    $-3/4 + 840i$ & $3.095 \cdots \times 10^{-59}$ \\
    $-\frac{1}{2} + 100i$ & $6.264 \cdots \times 10^{-138}$ \\
    $-\frac{1}{2} + 300i$ & $7.310 \cdots \times 10^{-133}$ \\
    $-\frac{1}{2} + 500i$ & $9.690 \cdots \times 10^{-100}$ \\
    $-\frac{1}{2} + 600i$ & $1.459 \cdots \times 10^{-85} $ \\
    $-\frac{1}{2} + 840i$ & $9.401 \cdots \times 10^{-60} $ \\
    $100i$ & $8.178 \cdots \times 10^{-139}$ \\
    $300i$ & $5.289 \cdots \times 10^{-134}$ \\
    \bottomrule
    \end{tabular}
    \end{subtable}
    &
    \begin{subtable}[t]{0.45\textwidth}
    \centering
    \begin{tabular}{c c}
    \toprule
    $s$ & \multicolumn{1}{c}{$|L(s, \chi_4) - F_{a_n,~ 702}(s)|$} \\
    \midrule
    $500i$ & $5.431 \cdots \times 10^{-101}$ \\
    $600i$ & $7.468 \cdots \times 10^{-87}$ \\
    $840i$ & $2.783 \cdots \times 10^{-57}$ \\
    $\frac{1}{2} + 300i$ & $3.827 \cdots \times 10^{-135}$ \\
    $\frac{1}{2} + 500i$ & $3.044 \cdots \times 10^{-102}$ \\
    $\frac{1}{2} + 600i$ & $8.867 \cdots \times 10^{-88}$ \\
    $\frac{1}{2} + 840i$ & $1.203 \cdots \times 10^{-57}$ \\
    $1 + 100i$ & $3.117 \cdots \times 10^{-140}$ \\
    $2 + 300i$ & $1.453 \cdots \times 10^{-138}$ \\
    $3 + 500i$ & $1.685 \cdots \times 10^{-108}$ \\
    $4 + 600i$ & $3.512 \cdots \times 10^{-97}$ \\
    $5 + 840i$ & $6.373 \cdots \times 10^{-70}$ \\
    \bottomrule
    \end{tabular}
    \end{subtable}
\end{tabular}
\end{table}

Indeed, the approximant \( F_{a_n,~ 702} \) offers a more accurate approximation of the zeros with imaginary parts in the range \([g_0, g_{699}]\) and beyond that, identifying at least the zeros \(\rho_1 = \frac{1}{2} + i \, 6.020948 \dots\) through \(\rho_{615} = \frac{1}{2} + i \, 825.2909367 \dots\), as detailed by Equation (80) to (88).

\begin{align}    
0 &= F_{a_n,~ 702}  \left( \rho_{1} +8.71181 \dots \times 10^{-141} -3.33865 \dots \times 10^{-140} i \right) \\
& \dots \nonumber \\ 
0 &= F_{a_n,~ 702}  \left( \rho_{100} +1.80565 \dots \times 10^{-141} -2.71973 \dots \times 10^{-140} i \right) \\
& \dots \nonumber \\ 
0 &= F_{a_n,~ 702}  \left( \rho_{200} +2.97736 \dots \times 10^{-137} +5.98309 \dots \times 10^{-136} i  \right) \\
& \dots \nonumber \\ 
0 &= F_{a_n,~ 702}  \left( \rho_{300} +3.00986 \dots \times 10^{-114} +8.69802 \dots \times 10^{-114} i  \right) \\
& \dots \nonumber \\ 
0 &= F_{a_n,~ 702}  \left( \rho_{400} -2.73477 \dots \times 10^{-95} +2.42604 \dots \times 10^{-95} i  \right) \\
& \dots \nonumber \\
0 &= F_{a_n,~ 702}  \left( \rho_{500} + 1.01226 \dots \times 10^{-77} +2.03857 \dots \times 10^{-78} i \right) \\
& \dots \nonumber \\
0 &= F_{a_n,~ 702}  \left( \rho_{615} + 7.49384 \dots \times 10^{-61} +3.50907 \dots \times 10^{-60} i \right) \\
& \dots \nonumber \\
0 &= F_{a_n,~ 702}  \left( \rho_{700} + 1.96222 \dots \times 10^{-48} +1.44749 \dots \times 10^{-48} i \right) \\
& \dots \nonumber \\
0 &= F_{a_n,~ 702}  \left( \rho_{800} + 3.93119 \dots \times 10^{-36} +3.02971 \dots \times 10^{-33} i \right) \\
& \dots \nonumber 
\end{align}

\subsubsection{Case of \( q = 3 \) and \( d = -3 \)}
For \( q = 3 \) and \( d = -3 \), the coefficients corresponding to indices that are multiples of 3 are set to zero. Taking \( k = 2 \), we set \( a_1 = \chi_3(1) = 1 \) and \( a_3 = \chi_3(2) = -1 \), and using \( M = 500 \) generalized Gram points (so that \( N = 502 \)), the coefficients of the approximant \( F_{a_n,~ 502} \) are depicted in Figure 2.3.2.1.

Table 1 provides a sample of the behavior of the approximant \( F_{a_n,~ 502} \), using the interpolating Gram points \( g_0 = 4.8301 \ldots \), \(\ldots\), and \( g_{499} = 658.4834 \ldots \).

\begin{figure}[ht]
 \renewcommand{\thefigure}{2.3.2.1}
  \centering
  \includegraphics[width=0.7\textwidth]{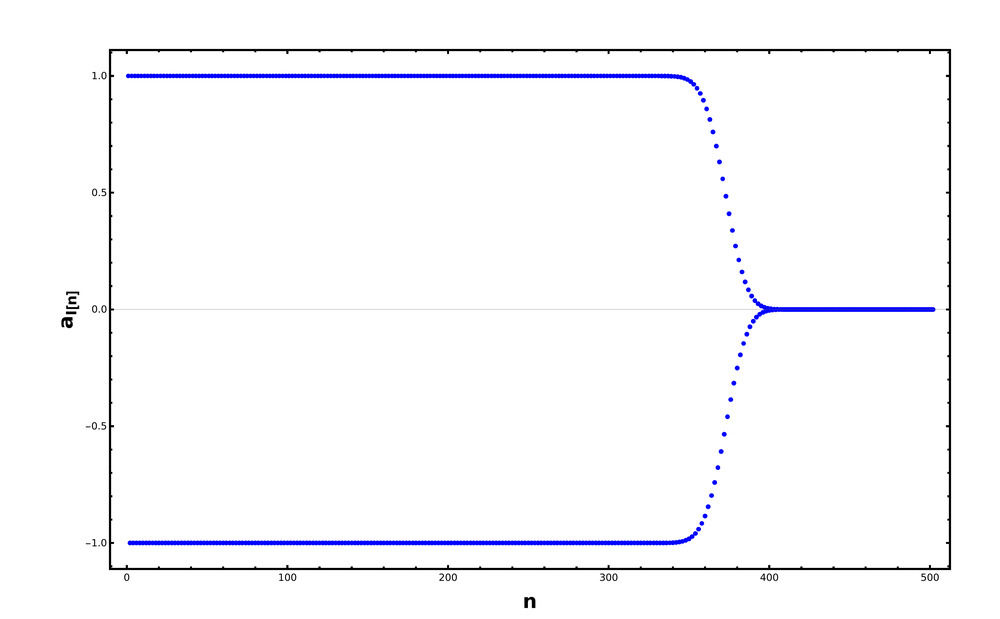}
  \caption{Coefficients \( a_{I[n]} \) of the approximant \( F_{a_n,~ 401}  \) when \( M=200 \), where \( I[n] \) is the \( n \)-th element of the list $I$.}
  \label{fig:thetaz7}
\end{figure}

The approximant \( F_{a_n,~ 502} \) provides a highly accurate approximation of the zeros of \( L(s,{\chi_3}) \), particularly within the range of imaginary parts from \( g_0 \) to \( g_{499} \), and even slightly beyond \( g_{499} \), as demonstrated in Equation (89) to (94). Remarkably, \( F_{a_n,~ 502} \) detects zeros on the critical line with imaginary parts between \( \Im(\rho_1) = 8.039737 \dots \) and \( \Im(\rho_{470}) = 660.877547  \dots \), and even beyond that.
\begin{table}
\centering
\caption{Approximation of $ L(s,{\chi_3})$ by $F_{a_n,~ 502}$.}

\begin{tabular}{cc}
    \begin{subtable}[t]{0.45\textwidth}
    \centering
    \begin{tabular}{c c}
    \toprule
    $s$ & \multicolumn{1}{c}{$|L(s,{\chi_3})- F_{a_n,~ 502} (s)|$} \\
    \midrule
    $-4 + 100i$ & $3.3439 \dots \times 10^{-89}$ \\
    $-3 + 300i$ & $7.1780 \dots \times 10^{-72}$ \\
    $-2 + 500i$ & $6.4511 \dots \times 10^{-48}$ \\
    $-1 + 640i$ & $8.1247 \dots \times 10^{-33}$ \\
    $-1/2 + 100i$ & $4.4959 \dots \times 10^{-84}$ \\
    $-1/2 + 300i$ & $2.9144 \dots \times 10^{-79}$ \\
    $-1/2 + 500i$ & $1.7477 \dots \times 10^{-47}$ \\
    $-1/2 + 640i$ & $4.6473 \dots \times 10^{-33}$ \\
    $100i$ & $2.4627 \dots \times 10^{-99}$ \\
    $300i$ & $2.4346 \dots \times 10^{-73}$ \\
    \bottomrule
    \end{tabular}
    \end{subtable}
    &
    \begin{subtable}[t]{0.45\textwidth}
    \centering
    \begin{tabular}{c c}
    \toprule
    $s$ & \multicolumn{1}{c}{$|L(s,{\chi_3})- F_{a_n,~ 502} (s)|$} \\
    \midrule
    $500i$ & $1.1310 \dots \times 10^{-48}$ \\
    $640i$ & $2.6583 \dots \times 10^{-34}$ \\
    $1/2 + 100i$ & $1.3536 \dots \times 10^{-86}$ \\
    $1/2 + 300i$ & $2.0341 \dots \times 10^{-80}$ \\
    $1/2 + 500i$ & $7.3200 \dots \times 10^{-50}$ \\
    $1/2 + 640i$ & $1.5207 \dots \times 10^{-35}$ \\
    $1 + 100i$ & $7.4648 \dots \times 10^{-88}$ \\
    $2 + 300i$ & $1.1872 \dots \times 10^{-77}$ \\
    $3 + 500i$ & $8.3190 \dots \times 10^{-56}$ \\
    $4 + 640i$ & $3.0517 \dots \times 10^{-44}$ \\
    \bottomrule
    \end{tabular}
    \end{subtable}
\end{tabular}
\end{table}

\begin{align}    
0 &= F_{a_n,~ 502}  \left( \rho_{1} + 2.16539 \dots \times 10^{-80} -1.17018 \dots \times 10^{-80} i \right) \\
& \dots \nonumber \\ 
0 &= F_{a_n,~ 502}  \left( \rho_{100} -1.91227 \dots \times 10^{-99} +7.38512 \dots \times 10^{-99}  i \right) \\
& \dots \nonumber \\ 
0 &= F_{a_n,~ 502}  \left( \rho_{200} +4.51875 \dots \times 10^{-73} +9.25327 \dots \times 10^{-74} i  \right) \\
& \dots \nonumber \\ 
0 &= F_{a_n,~ 502}  \left( \rho_{300} +3.22646 \dots \times 10^{-64} +1.67801 \dots \times 10^{-64} i  \right) \\
& \dots \nonumber \\ 
0 &= F_{a_n,~ 502}  \left( \rho_{400} -2.18916 \dots \times 10^{-43} +5.99374 \dots \times 10^{-43} i  \right) \\
& \dots \nonumber \\ 
0 &= F_{a_n,~ 502}  \left( \rho_{460} +  4.56864 \dots \times 10^{-36} +1.12057 \dots \times 10^{-35} i  \right) \\
& \dots \nonumber \\ 
0 &= F_{a_n,~ 502}  \left( \rho_{500} + 4.17372 \dots \times 10^{-34} +1.64367 \dots \times 10^{-34} i \right) \\
& \dots \nonumber
\end{align}

If we use more interpolating points, we can improve the accuracy of the approximation. For instance, using the Gram points \( g_0 = 4.8301 \ldots \), \(\ldots\), and \( g_{699} = 871.5669 \ldots \), Table 2 provides a sample of the behavior of the approximant \( F_{a_n,~ 702} \).

\begin{table}[ht]
\centering
\caption{Approximation of $L(s, \chi_3)$ by $F_{a_n,~ 702}$.}

\begin{tabular}{cc}
    \begin{subtable}[t]{0.45\textwidth}
    \centering
    \begin{tabular}{c c}
    \toprule
    $s$ & \multicolumn{1}{c}{$|L(s, \chi_3) - F_{a_n,~ 702}(s)|$} \\
    \midrule
    $-4 + 100i$ & $1.178 \cdots \times 10^{-131}$ \\
    $-3 + 300i$ & $3.695 \cdots \times 10^{-127}$ \\
    $-2 + 500i$ & $5.508 \cdots \times 10^{-96}$ \\
    $-1 + 600i$ & $2.853 \cdots \times 10^{-84}$ \\
    $-3/4 + 840i$ & $3.095 \cdots \times 10^{-59}$ \\
    $-\frac{1}{2} + 100i$ & $6.264 \cdots \times 10^{-138}$ \\
    $-\frac{1}{2} + 300i$ & $7.310 \cdots \times 10^{-133}$ \\
    $-\frac{1}{2} + 500i$ & $9.690 \cdots \times 10^{-100}$ \\
    $-\frac{1}{2} + 600i$ & $1.459 \cdots \times 10^{-85} $ \\
    $-\frac{1}{2} + 840i$ & $9.401 \cdots \times 10^{-60} $ \\
    $100i$ & $8.178 \cdots \times 10^{-139}$ \\
    $300i$ & $5.289 \cdots \times 10^{-134}$ \\
    \bottomrule
    \end{tabular}
    \end{subtable}
    &
    \begin{subtable}[t]{0.45\textwidth}
    \centering
    \begin{tabular}{c c}
    \toprule
    $s$ & \multicolumn{1}{c}{$|L(s, \chi_3) - F_{a_n,~ 702}(s)|$} \\
    \midrule
    $500i$ & $5.431 \cdots \times 10^{-101}$ \\
    $600i$ & $7.468 \cdots \times 10^{-87}$ \\
    $840i$ & $2.783 \cdots \times 10^{-57}$ \\
    $\frac{1}{2} + 300i$ & $3.827 \cdots \times 10^{-135}$ \\
    $\frac{1}{2} + 500i$ & $3.044 \cdots \times 10^{-102}$ \\
    $\frac{1}{2} + 600i$ & $8.867 \cdots \times 10^{-88}$ \\
    $\frac{1}{2} + 840i$ & $1.203 \cdots \times 10^{-57}$ \\
    $1 + 100i$ & $3.117 \cdots \times 10^{-140}$ \\
    $2 + 300i$ & $1.453 \cdots \times 10^{-138}$ \\
    $3 + 500i$ & $1.685 \cdots \times 10^{-108}$ \\
    $4 + 600i$ & $3.512 \cdots \times 10^{-97}$ \\
    $5 + 840i$ & $6.373 \cdots \times 10^{-70}$ \\
    \bottomrule
    \end{tabular}
    \end{subtable}
\end{tabular}
\end{table}

Indeed, the approximant \( F_{a_n,~ 702} \) offers a more accurate approximation of $L(s, \chi_3)$ and its zeros with imaginary parts in the range \([g_0, g_{699}]\) and beyond that, identifying at least the zeros \(\rho_1 = \frac{1}{2} + i \, 8.039737 \dots\) through \(\rho_{636} = \frac{1}{2} + i \, 870.903928  \dots\), as detailed by Table 4 and Equation (96) to (104).

\begin{align}     
0 &= F_{a_n,~ 702}  \left( \rho_{1} -7.7679 \dots \times 10^{-130} -1.66431 \dots \times 10^{-129} i \right) \\
& \dots \nonumber \\ 
0 &= F_{a_n,~ 702}  \left( \rho_{100} -2.63906 \dots \times 10^{-145} +3.64137 \dots \times 10^{-146} i \right) \\
& \dots \nonumber \\ 
0 &= F_{a_n,~ 702}  \left( \rho_{200} -5.43641 \dots \times 10^{-129} -1.11324 \dots \times 10^{-129} i  \right) \\
& \dots \nonumber \\ 
0 &= F_{a_n,~ 702}  \left( \rho_{300} -6.18454 \dots \times 10^{-108} -3.21644 \dots \times 10^{-108} i  \right) \\
& \dots \nonumber \\ 
0 &= F_{a_n,~ 702}  \left( \rho_{400} +3.85683 \dots \times 10^{-99} -1.05597 \dots \times 10^{-98} i  \right) \\
& \dots \nonumber \\ 
0 &= F_{a_n,~ 702}  \left( \rho_{500} -1.15064 \dots \times 10^{-74} -4.53138 \dots \times 10^{-74} i \right) \\
& \dots \nonumber \\
0 &= F_{a_n,~ 702}  \left( \rho_{636} +1.76602 \dots \times 10^{-54} -3.86543 \dots \times 10^{-54} i \right) \\
& \dots \nonumber 
\end{align}
\begin{align}
0 &= F_{a_n,~ 702}  \left( \rho_{700} -3.42705 \dots \times 10^{-49} -2.21712 \dots \times 10^{-49} i \right) \\
& \dots \nonumber \\  
0 &= F_{a_n,~ 702}  \left( \rho_{800} -3.63141 \dots \times 10^{-32} +1.59907 \dots \times 10^{-32} i \right) \\
& \dots \nonumber 
\end{align}

\section{Future Directions} 
\begin{itemize}
    \item \textbf{Discovering new $L$-functions:} From the results of this paper, one observes that the solutions to $\theta(t, \chi)$ totally determine the L-function corresponding to the function $\theta$ (e.g. the coefficients in figure 2.2.1.2 determine all the coefficients by periodicity). Thus, this presents an interesting question: given a list of real numbers that are the gram points of some mysterious $\theta$- function corresponding to some mysterious $L$-function with coefficients satisfying some periodicity condition, can we discover $L$ by applying an experiment similar to the ones carried out in this paper?
    \item \textbf{Employing more advanced feature selection methods from machine learning to make the approach of this paper work for Dirichlet $L$- functions with very large modulus}
\end{itemize}

\section*{Acknowledgment} The author would like to express gratitude to Professor Ghaith Hiary for his assistance in answering some queries during the course of this research. The author received no financial support for the research and authorship of this article.


\begin{thebibliography}{99} 

\bibitem{DaviesWilkes1965}
H. Davies and M. V. Wilkes,
\emph{An approximate functional equation for Dirichlet L-functions},
\emph{Proceedings of the Royal Society of London A}, vol. 284, pp. 224-236, 1965.

\bibitem{Spira1969}
R. Spira,
\emph{Calculation of Dirichlet L-Functions},
\emph{Mathematics of Computation}, vol. 23, no. 107, pp. 489-513, 1969. 
Available: JSTOR, \url{https://doi.org/10.2307/2004376}. 
Accessed: Aug. 13, 2024.

\bibitem{Rumely1993}
R. Rumely,
\emph{Supplement to Numerical Computations Concerning the ERH},
\emph{Mathematics of Computation}, vol. 61, no. 203, pp. S17-23, 1993.
Available: JSTOR, \url{https://doi.org/10.2307/2152974}.
Accessed: Aug. 13, 2024.

\bibitem{Lavrik1968}
A. F. Lavrik,
\emph{Functional and approximate functional equations of the Dirichlet function},
\emph{Mathematical Notes of the Academy of Sciences of the USSR}, vol. 3, no. 5, pp. 388-393, May 1968.
Available: \url{https://doi.org/10.1007/BF01150995}.

\bibitem{Riemann1859}
B. Riemann,
\emph{Ueber die Anzahl der Primzahlen unter einer gegebenen Grösse},
\emph{Monatsberichte der Königlichen Preussische Akademie der Wissenschaften zu Berlin}, pp. 671-680, 1859.

\bibitem{HiaryOdlyzko2012}
G. A. Hiary and A. M. Odlyzko,
\emph{The zeta function on the critical line: Numerical evidence for moments and random matrix theory models},
\emph{Mathematics of Computation}, vol. 81, no. 279, pp. 1723–1752, 2012.

\bibitem{Gram1903}
J.-P. Gram,
\emph{Note sur les zéros de la fonction $\zeta(s)$ de Riemann},
\emph{Acta Mathematica}, vol. 27, no. 1, pp. 289-304, 1903.
\url{https://doi.org/10.1007/BF02421310}

\bibitem{Hutchinson1925}
J. I. Hutchinson,
\emph{On the Roots of the Riemann Zeta-Function},
\emph{Transactions of the American Mathematical Society}, vol. 27, no. 1, pp. 49-60, 1925.

\bibitem{Titchmarsh1935}
E. C. Titchmarsh,
\emph{The zeros of the Riemann zeta-function},
\emph{Proceedings of the Royal Society of London}, vol. 151, pp. 234-255, 1935; vol. 157, pp. 261-263, 1936.


\bibitem{Lehmer1956}
D. H. Lehmer,
\emph{On the roots of the Riemann zeta-function},
\emph{Acta Mathematica}, vol. 95, pp. 291-298, 1956.



\bibitem{MATIYASEVICH2020460}
Yu. Matiyasevich,
\emph{Plausible ways for calculating the Riemann zeta function via the Riemann–Siegel theta function},
\emph{Journal of Number Theory}, vol. 207, pp. 460-471, 2020.
\url{https://doi.org/10.1016/j.jnt.2019.07.021}

\bibitem{matiyasevich2013}
Yuri Matiyasevich,
\emph{Calculation of Riemann’s zeta function via interpolating determinants},
2013, Preprint 2013-18, Max Planck Institute for Mathematics in Bonn,
Available online: \url{http://www.mpim-bonn.mpg.de/preblob/5368}.

\bibitem{Trefethen2021}
Lloyd N. Trefethen.
\textit{Approximation Theory and Approximation Practice, Extended Edition}.
SIAM, 2021.
\bibitem{Ivic1985}
Aleksandar Ivic,
\textit{The Riemann Zeta-Function: Theory and Applications}.
Dover Publications, 1985.
\bibitem{Edwards1974}
H. M. Edwards,
\textit{Riemann's Zeta Function},
New York: Dover Publications, 1974. ISBN 978-0-486-41740-0. MR 0466039.
\bibitem{Gabcke1979}
W. Gabcke,
\textit{Neue Herleitung und explizierte Restabschätzung der Riemann-Siegel-Formel}.
Thesis, University of Göttingen, 1979. Revised version (eDiss Göttingen 2015).
\bibitem{Davenport2000}
H. Davenport,
\textit{Multiplicative Number Theory},
3rd edition, Graduate Texts in Mathematics, vol. 74, Springer, 2000, edited by Hans Heilbronn, ISBN 978-0-387-95097-5.
\bibitem{Apostol1976}
T. M. Apostol,
\textit{Introduction to Analytic Number Theory}.
Springer-Verlag, New York, 1976.
\bibitem{MontgomeryVaughan2006}
H. L. Montgomery and R. C. Vaughan,
\textit{Multiplicative Number Theory. I. Classical Theory}.
Cambridge Tracts in Advanced Mathematics, vol. 97, Cambridge University Press, 2006.
ISBN 978-0-521-84903-6.
\bibitem{Kronecker1885}
L. Kronecker,
\textit{Zur Theorie der elliptischen Funktionen},
Sitzungsberichte der Königlich Preussischen Akademie der Wissenschaften zu Berlin, 
pp. 761--784, 1885.
\bibitem{Cramer1750}
G. Cramer, \textit{Introduction à l'Analyse des lignes Courbes algébriques}, Geneva: Europeana, 1750, pp. 656–659. Retrieved 2012-05-18.
\bibitem{IbukiyamaKaneko2014}
Ibukiyama, T., Kaneko, M. (2014). The Euler–Maclaurin Summation Formula and the Riemann Zeta Function. In: Bernoulli Numbers and Zeta Functions. Springer Monographs in Mathematics. Springer, Tokyo.
\bibitem{Gourdon2003NumericalEO}
Gourdon, X., Sebah, P. (2003). Numerical evaluation of the Riemann Zeta-function. \url{https://api.semanticscholar.org/CorpusID:125658246}.
\bibitem{Odlyzko1988Fast}
Odlyzko, A. M., Schönhage, A. (1988). Fast algorithms for multiple evaluations of the Riemann zeta function. Trans. Amer. Math. Soc. 309 (2): 797–809. \url{https://doi.org/10.2307/2000939}. JSTOR 2000939. MR 0961614.
\bibitem{PlattTrudgian2021}
Platt, D., Trudgian, T. (2021). The Riemann hypothesis is true up to \(3 \cdot 10^{12}\). Bulletin of the London Mathematical Society, 53. \url{https://doi.org/10.1112/blms.12460}.
\bibitem{DeBoor2001}
C. De Boor, \emph{A Practical Guide to Splines}, Springer-Verlag, 2001. ISBN 978-0-387953-66-9.
\bibitem{Saad1986}
Y. Saad and M. H. Schultz, \emph{GMRES: A Generalized Minimal Residual Algorithm for Solving Nonsymmetric Linear Systems}, SIAM J. Sci. Stat. Comput., vol. 7, no. 3, pp. 856--869, 1986. DOI: 10.1137/0907058.
\bibitem{Saad1996}
Y. Saad, \emph{Iterative Methods for Sparse Linear Systems}, 1st ed., PWS, Boston, 1996. ISBN 978-0-534-94776-7.
\bibitem{Nifa2017}
Naoufal Nifa,
\textit{Solveurs performants pour l'optimisation sous contraintes en identification de paramètres},
PhD thesis, 2017.
\bibitem{hestenes1952}
Magnus R. Hestenes and Eduard Stiefel, 
\emph{Methods of conjugate gradients for solving linear systems}, 
1952, Journal of Research of the National Bureau of Standards, Vol. 49, pp. 409--435, 
Available online: \url{https://api.semanticscholar.org/CorpusID:2207234}.
\bibitem{fonti2017}
Fonti, Valeria, and Eduard Belitser. "Feature selection using lasso." VU Amsterdam Research Paper in Business Analytics 30 (2017): 1-25.



\end{thebibliography}
\end{document}